\documentclass[11pt,twoside]{amsart}
\usepackage{amsmath, amsthm, amscd, amsfonts, amssymb, graphicx}
\usepackage[bookmarksnumbered, plainpages]{hyperref}
\usepackage{mathrsfs}
\UseRawInputEncoding
\newtheorem{thm}{Theorem}[section]

\newtheorem{lem}[thm]{Lemma}

\newtheorem*{Proof}{Proof}

\numberwithin{equation}{section}

\begin{document}

\title{\bf Conformal Ricci collineations associated to the Levi-Civita connection on three-dimensional Lorentzian Lie groups}
\author{Yanli Wang}

\thanks{{\scriptsize
\hskip -0.4 true cm \textit{2010 Mathematics Subject Classification:}
53C40; 53C42.
 }}

\maketitle

\begin{abstract}
 In this paper, we determine all conformal Ricci collineations associated to the Levi-Civita connection on three-dimensional Lorentzian Lie groups .
\end{abstract}

\vskip 0.2 true cm


\pagestyle{myheadings}
\markboth{\rightline {\scriptsize Wang}}
         {\leftline{\scriptsize Conformal Ricci collineations}}

\bigskip
\bigskip

\section{ Introduction}
\indent Symmetry can be regarded as a one-parameter group of diffeomorphism of spacetime, and it preserves mathematical and physical quantity. Symmetric systems are relatively simple in form and have special properties. Symmetry plays an irreplaceable role in spacetime. In 1969, based on the different geometric objects' invariance properties, Katzin defined them as those vector fields X, such that leave the various relevant geometric quantities and classified Ricci collineations, curvature collineations in \cite{Ka}. Collineations are symmetry properties of space-times, so Ricci collineations i.e. Ricci symmetry, which is defined by $(L_{\xi}Ric)=0$. Because Ricci Collineations' close connection with the energy-momentum tensor, inspired more mathematicians' research interesting. In \cite{Du}, Duggal generalized Ricci collineation to a new symmetry called conformal Ricci collineation. W. K\"{u}hnel defined a more general conformal Ricci collineation by the combination of $(L_{\xi}g)=2\sigma g$ and $(L_{\xi}Ric)=2\lambda g$ in \cite{Ku}. When $\lambda=0$, conformal Ricci collineation is reduced to Ricci collineation. Ricci Collineations and conformal Ricci collineations have been discussed and determined in more different spacetimes and for various other ricci tensor in \cite{CN, AC, CU,MU}.\\
\indent In \cite{Ca1,CP}, Calvaruso and Cordero divided three-dimensional Lorentzian Lie groups into $\{G_i\}_{i=1,\cdots,7}$, and geted corresponding Lie algebra $\{\mathfrak{g}\}_{i=1,\cdots,7.}$. By this classification, Batat studied algebraic Ricci solitons on three-dimensional Lorentzian Lie groups in \cite{W}. In this paper, similar to \cite{W}, we compute the Ricci tensor of $\{G_i\}_{i=1,\cdots,7}$ associated to the Levi-Civita connection $\nabla^{L}$, and then determine all conformal Ricci collineations associated to the Levi-Civita connection on $\{G_i\}_{i=1,\cdots,7}$ . Where $\{G_i\}_{i=1,\cdots,4}$ are three-dimensional unimodular Lorentzian Lie groups; $\{G_i\}_{i=5,\cdots,7}$ are three-dimensional non-unimodular Lorentzian Lie groups.\\
\indent In section 2, we first recall the definition of the Ricci collineation, and generalize to conformal Ricci collineation, then we determine conformal Ricci collineations associated to the Levi-Civita connection on three-dimensional Lorentzian Unimodular Lie groups. In section 3, we determine conformal Ricci collineations associated to the Levi-Civita connection on three-dimensional Lorentzian Unimodular Lie groups.


\vskip 1 true cm

\section{ Conformal Ricci collineations associated to the Levi-Civita connection on three-dimensional Lorentzian Non-unimodular Lie groups}

\indent Calvaruso and Cordero classified three-dimensional Lorentzian Lie groups into $\{G_i\}_{i=1,\cdots,7}$ equipped with a left-invariant Lorentzian metric $g$ in \cite{Ca1,CP}(see Theorem 2.1 and Theorem 2.2 in \cite{W}), and geted corresponding Lie algebra $\{\mathfrak{g}\}_{i=1,\cdots,7.}$. In \cite{W}, Batat studied algebraic Ricci solitons of three-dimensional Lorentzian Lie groups, similar to his calculations, we compute the Ricci tensor of $\{G_i\}_{i=1,\cdots,7}$ associated to the Levi-Civita connection $\nabla^{L}$. We define the curvature tensor of the Levi-Civita connection $\nabla^{L}$:
\begin{equation}
R^{L}(X,Y)Z=\nabla^{L}_X\nabla^{L}_YZ-\nabla^{L}_Y\nabla^{L}_XZ-\nabla^{L}_{[X,Y]}Z.
\end{equation}
The Ricci tensor of $(G_i,g)$ associated to the Levi-Civita connection $\nabla^{L}$ is defined by
\begin{equation}
\rho^{L}(X,Y)=-g(R^{L}(X,e_1)Y,e_1)-g(R^{L}(X,e_2)Y,e_2)+g(R^{L}(X,e_3)Y,e_3).
\end{equation}
where $g(e_i,e_j)=0,i\neq j;~g(e_i,e_i)=1,i=1,2;~g(e_3,e_3)=-1.$ Let
\begin{equation}
Ric^{L}(X,Y)=\frac{\rho^{L}(X,Y)+\rho^{L}(Y,X)}{2}.
\end{equation}
we get the symmetric Ricci tensor.\\
\indent We know that for vector $X,Y,V$, the Lie derivative of the Ricci tensor $Ric^{L}$ associated to $V$ can be defined:
\begin{equation}
(L_ZRic^{L})(X,Y):=Z[Ric^{L}(X,Y)]-Ric^{L}([Z,X],Y)-Ric^{L}(X,[Z,Y])
\end{equation}
$(G_i,g)$ admits left-invariant Ricci collineations associated to the Levi-Civita connection $\nabla^{L}$ if and only if it satisfies
$(L_VRic^{L})(X,Y)=0$.
Further, conformal Ricci collineations can be discribed:
\begin{thm} $(G_i,g)$ admits conformal Ricci collineations associated to the Levi-Civita connection $\nabla^{L}$ if and only if it satisfies
\begin{equation}
(L_VRic^{L})=2\lambda g,
\end{equation}
where $V=\lambda_1e_1+\lambda_2e_2+\lambda_3e_3$ is a left-invariant vector field and $\lambda_1,\lambda_2,\lambda_3$ are real numbers.
\end{thm}

\vskip 0.5 true cm
\noindent{\bf 2.1 Conformal Ricci collineations of $G_1$ }\\
\vskip 0.5 true cm
 By \cite{W}, we have the following Lie algebra of $G_1$ satisfies
\begin{equation}
[e_1,e_2]=\alpha e_1-\beta e_3,~~[e_1,e_3]=-\alpha e_1-\beta e_2,~~[e_2,e_3]=\beta e_1+\alpha e_2+\alpha e_3,~~\alpha\neq 0.
\end{equation}
where $e_1,e_2,e_3$ is a pseudo-orthonormal basis, with $e_3$ timelike.
\begin{lem}The Ricci tensor of  $(G_1,g)$ associated to the Levi-Civita connection $\nabla^{L}$ is determined by
\begin{align}
&Ric^{L}(e_1,e_1)=-\frac{\beta^2}{2},~~~Ric^{L}(e_1,e_2)=-\alpha\beta,~~~Ric^{L}(e_1,e_3)=\alpha\beta,\\\notag
&Ric^{L}(e_2,e_2)=-2\alpha^2-\frac{\beta^2}{2},~~~Ric^{L}(e_2,e_3)=2\alpha^2,~~~Ric^{L}(e_3,e_3)=-2\alpha^2+\frac{\beta^2}{2}.
\end{align}
\end{lem}
 Let $V=\lambda_1e_1+\lambda_2e_2+\lambda_3e_3$ is a left-invariant vector field, we can get the Lie derivative of the Ricci tensor $Ric^{L}$ associated to $V$.
\begin{lem}
\begin{align}
&(L_VRic^{L})(e_1,e_1)=-3\alpha\beta^2\lambda_2+3\alpha\beta^2\lambda_3,\\\notag
&(L_VRic^{L})(e_1,e_2)=\frac{3}{2}\alpha\beta^2\lambda_1-3\alpha^2\beta\lambda_2+3\alpha^2\beta\lambda_3,\\\notag
&(L_VRic^{L})(e_1,e_3)=-\frac{3}{2}\alpha\beta^2\lambda_1+3\alpha^2\beta\lambda_2+\alpha^2\beta\lambda_3,\\\notag
&(L_VRic^{L})(e_2,e_2)=6\alpha^2\beta\lambda_1-3\alpha\beta^2\lambda_3,\\\notag
&(L_VRic^{L})(e_2,e_3)=-6\alpha^2\beta\lambda_1+\dfrac{3}{2}\alpha\beta^2\lambda_2+\dfrac{3}{2}\alpha\beta^2\lambda_3,\\\notag
&(L_VRic^{L})(e_3,e_3)=6\alpha^2\beta\lambda_1-3\alpha\beta^2\lambda_2.
\end{align}
\end{lem}
Then, if $V$ is a conformal Ricci collineation associated to the Levi-Civita connection, by Lemma 2.3 and Theorem 2.1, we have the following equations:
\begin{eqnarray}
       \begin{cases}
       -3\alpha\beta^2\lambda_2+3\alpha\beta^2\lambda_3=2\lambda\\[2pt]
       \dfrac{3}{2}\alpha\beta^2\lambda_1-3\alpha^2\beta\lambda_2+3\alpha^2\beta\lambda_3=0\\[2pt]
       -\dfrac{3}{2}\alpha\beta^2\lambda_1+3\alpha^2\beta\lambda_2+\alpha^2\beta\lambda_3=0\\[2pt]
       6\alpha^2\beta\lambda_1-3\alpha\beta^2\lambda_3=2\lambda\\[2pt]
       -6\alpha^2\beta\lambda_1+\dfrac{3}{2}\alpha\beta^2\lambda_2+\dfrac{3}{2}\alpha\beta^2\lambda_3=0\\[2pt]
       6\alpha^2\beta\lambda_1-3\alpha\beta^2\lambda_2=-2\lambda\\[2pt]
       \end{cases}
\end{eqnarray}
By solving (2.9) , we get
\begin{thm}
$(G_1, g, V)$ admits conformal Ricci collineations associated to the Levi-Civita connection if and only if:
 $\alpha\neq0,\beta= 0,\lambda= 0,~\mathscr{V}_{\mathscr{R}{C}}=<e_1,e_2,e_3>$, where $\mathscr{V}_{\mathscr{R}{C}}$ is the vector space of conformal Ricci collineations on $(G_1, g, V)$.
\end{thm}
\begin{Proof}
We know that $\alpha\neq 0$. By the second and third equations, we get $\alpha^2\beta\lambda_3=0$, then $\beta\lambda_3=0$, so\\
case {\rm 1)} If $\lambda_3=0$, by {\rm (2.9)},
\begin{eqnarray}
       \begin{cases}
       -3\alpha\beta^2\lambda_2=2\lambda\\[2pt]
       \dfrac{3}{2}\alpha\beta^2\lambda_1-3\alpha^2\beta\lambda_2=0\\[2pt]
       6\alpha^2\beta\lambda_1=2\lambda\\[2pt]
       -6\alpha^2\beta\lambda_1+\dfrac{3}{2}\alpha\beta^2\lambda_2=0\\[2pt]
       6\alpha^2\beta\lambda_1-3\alpha\beta^2\lambda_2=-2\lambda\\[2pt]
       \end{cases}
\end{eqnarray}
then we get $\lambda=\beta=0$.\\
case {\rm 2)} If $\lambda_3\neq 0$, then we get $\beta=0$. {\rm (2.9)} holds.\\
So we have  $\alpha\neq0,\beta= 0,\lambda= 0$.
\end{Proof}

\vskip 0.5 true cm
\noindent{\bf 2.2 Conformal Ricci collineations of $G_2$ }\\
\vskip 0.5 true cm
 By \cite{W}, we have the following Lie algebra of $G_2$ satisfies
\begin{equation}
[e_1,e_2]=\gamma e_2-\beta e_3,~~[e_1,e_3]=-\beta e_2-\gamma e_3,~~[e_2,e_3]=\alpha e_1,~~\gamma\neq 0.
\end{equation}
where $e_1,e_2,e_3$ is a pseudo-orthonormal basis, with $e_3$ timelike.
\begin{lem}The Ricci tensor of  $(G_2,g)$ associated to the Levi-Civita connection $\nabla^{L}$ is determined by
\begin{align}
&Ric^{L}(e_1,e_1)=-\frac{\alpha^2}{2}-2\gamma^2,~~~Ric^{L}(e_1,e_2)=0,~~~Ric^{L}(e_1,e_3)=0,\\\notag
&Ric^{L}(e_2,e_2)=\frac{\alpha^2}{2}-\alpha\beta,~~~Ric^{L}(e_2,e_3)=-\alpha\gamma+2\beta\gamma,~~~Ric^{L}(e_3,e_3)=-\frac{\alpha^2}{2}+\alpha\beta.
\end{align}
\end{lem}
Similar to the previous lemma, we can get the Lie derivative of the Ricci tensor $Ric^{L}$ associated to $V$.
\begin{lem}
\begin{align}
&(L_VRic^{L})(e_1,e_1)=0,\\\notag
&(L_VRic^{L})(e_1,e_2)=(\dfrac{1}{2}\alpha^2\gamma-2\beta^2\gamma)\lambda_2+(-\dfrac{1}{2}\alpha^3-\alpha\gamma^2-2\beta\gamma^2+\alpha\beta^2-\dfrac{1}{2}\alpha^2\beta)\lambda_3,\\\notag
&(L_VRic^{L})(e_1,e_3)=(\dfrac{1}{2}\alpha^3+\alpha\gamma^2+2\beta\gamma^2-\alpha\beta^2+\dfrac{1}{2}\alpha^2\beta)\lambda_2+(\dfrac{1}{2}\alpha^2\gamma-2\beta^2\gamma)\lambda_3,\\\notag
&(L_VRic^{L})(e_2,e_2)=(-\alpha^2\gamma+4\beta^2\gamma)\lambda_1,\\\notag
&(L_VRic^{L})(e_2,e_3)=0,\\\notag
&(L_VRic^{L})(e_3,e_3)=(-\alpha^2\gamma+4\beta^2\gamma)\lambda_1.
\end{align}
\end{lem}
Then, if $V$ is a conformal Ricci collineation associated to the Levi-Civita connection, by Lemma 2.6 and Theorem 2.1, we have the following equations:
\begin{eqnarray}
       \begin{cases}
       0=2\lambda\\[2pt]
      (\dfrac{1}{2}\alpha^2\gamma-2\beta^2\gamma)\lambda_2+(-\dfrac{1}{2}\alpha^3-\alpha\gamma^2-2\beta\gamma^2+\alpha\beta^2-\dfrac{1}{2}\alpha^2\beta)\lambda_3=0\\[2pt]
      (\dfrac{1}{2}\alpha^3+\alpha\gamma^2+2\beta\gamma^2-\alpha\beta^2+\dfrac{1}{2}\alpha^2\beta)\lambda_2+(\dfrac{1}{2}\alpha^2\gamma-2\beta^2\gamma)\lambda_3=0\\[2pt]
      (-\alpha^2\gamma+4\beta^2\gamma)\lambda_1=2\lambda\\[2pt]
\end{cases}
\end{eqnarray}
By solving (2.14) , we get
\begin{thm}
$(G_2, g, V)$ admits conformal Ricci collineations associated to the Levi-Civita connection if and only if one of the following holds:\\
{\rm (1)} $\lambda=0,\alpha=2\beta,\gamma\neq 0, \gamma\neq 0,\beta\neq 0, ~\mathscr{V}_{\mathscr{R}{C}}=<e_1>$,\\
{\rm (2)} $\lambda=0,\alpha=-2\beta,\gamma\neq 0, ~\mathscr{V}_{\mathscr{R}{C}}=<e_1,e_2,e_3>$,\\
 where $\mathscr{V}_{\mathscr{R}{C}}$ is the vector space of conformal Ricci collineations on $(G_2, g, V)$.
\end{thm}
\begin{Proof}
We know that $\gamma\neq 0,\lambda=0$, \\
case {\rm 1)} If $\dfrac{1}{2}\alpha^2\gamma-2\beta^2\gamma=0$, i.e. $\alpha^2=4\beta^2$, then $(-\dfrac{1}{2}\alpha^3-\alpha\gamma^2-2\beta\gamma^2+\alpha\beta^2-\dfrac{1}{2}\alpha^2\beta)\lambda_3=0$. \\
case {\rm 1-1)} If $\alpha=2\beta$, then we have $\beta(\beta^2+\gamma^2)\lambda_3=0$, i.e. $\beta\lambda_3=0$. When $\beta\neq 0$, we have $\lambda_2=\lambda_3=0$, we get {\rm (1)}. When $\beta= 0$, we have $\alpha=0$, {\rm (2.14)} holds. This situation falls into {\rm (2)}.\\
case {\rm 1-2)} If $\alpha=-2\beta$, {\rm (2.14)} holds. This situation falls into {\rm (2)}.\\
case {\rm 2)} If $\dfrac{1}{2}\alpha^2\gamma-2\beta^2\gamma\neq 0$, i.e. $\alpha^2\neq 4\beta^2$, then $\lambda_2=\dfrac{b}{a}\lambda_3$, by the third equation we have
\begin{align}\notag
&b\dfrac{b}{a}\lambda_3+a\lambda_3=0,\notag
\end{align}
 i.e. $(b^2+a^2)\lambda_3=0$, we have $a=b=0$, where $b=\dfrac{1}{2}\alpha^3+\alpha\gamma^2+2\beta\gamma^2-\alpha\beta^2+\dfrac{1}{2}\alpha^2\beta, a=\dfrac{1}{2}\alpha^2\gamma-2\beta^2\gamma$. This is a contradiction.
\end{Proof}

\vskip 0.5 true cm
\noindent{\bf 2.3 Conformal Ricci collineations of $G_3$ }\\
\vskip 0.5 true cm
 By \cite{W}, we have the following Lie algebra of $G_3$ satisfies
\begin{equation}
[e_1,e_2]=-\gamma e_3,~~[e_1,e_3]=-\beta e_2,~~[e_2,e_3]=\alpha e_1.
\end{equation}
where $e_1,e_2,e_3$ is a pseudo-orthonormal basis, with $e_3$ timelike.
\begin{lem}The Ricci tensor of  $(G_3,g)$ associated to the Levi-Civita connection $\nabla^{L}$ is determined by
\begin{align}
&Ric^{L}(e_1,e_1)=-\alpha a_1-\beta a_2-\gamma a_3,~~~Ric^{L}(e_1,e_2)=0,~~~Ric^{L}(e_1,e_3)=0,\\\notag
&Ric^{L}(e_2,e_2)=\alpha a_1+\beta a_2-\gamma a_3 ,~~~Ric^{L}(e_2,e_3)=0,~~~Ric^{L}(e_3,e_3)=-\alpha a_1+\beta a_2-\gamma a_3.
\end{align}
\end{lem}
Similar to the previous lemma, we can get the Lie derivative of the Ricci tensor $Ric^{L}$ associated to $V$.
\begin{lem}
\begin{align}
&(L_VRic^{L})(e_1,e_1)=0,\\\notag
&(L_VRic^{L})(e_1,e_2)=(\alpha-\beta)[\gamma^2-(\alpha+\beta)^2]\lambda_3,\\\notag
&(L_VRic^{L})(e_1,e_3)=(\alpha-\gamma)[\beta^2-(\alpha+\gamma)^2]\lambda_2,\\\notag
&(L_VRic^{L})(e_2,e_2)=0,\\\notag
&(L_VRic^{L})(e_2,e_3)=(\beta-\gamma)[\alpha^2-(\beta+\gamma)^2]\lambda_1,\\\notag
&(L_VRic^{L})(e_3,e_3)=0.
\end{align}
\end{lem}
Then, if $V$ is a conformal Ricci collineation associated to the Levi-Civita connection, by Lemma 2.9 and Theorem 2.1, we have the following equations:
\begin{eqnarray}
       \begin{cases}
       0=2\lambda\\[2pt]
       (\alpha-\beta)[\gamma^2-(\alpha+\beta)^2]\lambda_3=0\\[2pt]
       (\alpha-\gamma)[\beta^2-(\alpha+\gamma)^2]\lambda_2=0\\[2pt]
       (\beta-\gamma)[\alpha^2-(\beta+\gamma)^2]\lambda_1=0\\[2pt]
       \end{cases}
\end{eqnarray}
By solving (2.18) , we get
\begin{thm}
$(G_3, g, V)$ admits conformal Ricci collineations associated to the Levi-Civita connection if and only if one of the following holds:\\
{\rm (1)} $\alpha+\beta+\gamma= 0,\lambda=0, \mathscr{V}_{\mathscr{R}{C}}=<e_1,e_2,e_3>$,\\
{\rm (2)} $\alpha+\beta+\gamma\neq 0,\alpha=\beta=\gamma\neq 0,\lambda=0, \mathscr{V}_{\mathscr{R}{C}}=<e_1,e_2,e_3>$,\\
{\rm (3)} $\alpha= 0,\beta=\gamma\neq 0,\lambda=0, \mathscr{V}_{\mathscr{R}{C}}=<e_1,e_2,e_3>$,\\
{\rm (4)} $\beta= 0,\alpha=\gamma\neq 0,\lambda=0, \mathscr{V}_{\mathscr{R}{C}}=<e_1,e_2,e_3>$,\\
{\rm (5)} $\gamma= 0,\alpha=\beta\neq 0,\lambda=0, \mathscr{V}_{\mathscr{R}{C}}=<e_1,e_2,e_3>$,\\
{\rm (6)} $\alpha+\beta+\gamma\neq 0,\alpha\neq 0,\beta=\gamma,\alpha\neq\gamma, \lambda=0,\mathscr{V}_{\mathscr{R}{C}}=<e_1>$,\\
{\rm (7)} $\alpha+\beta+\gamma\neq 0,\beta\neq 0,\alpha=\gamma,\alpha\neq\beta, \lambda=0,\mathscr{V}_{\mathscr{R}{C}}=<e_2>$,\\
{\rm (8)} $\alpha+\beta+\gamma\neq 0,\gamma\neq 0,\alpha=\beta,\beta\neq\gamma,\lambda=0, \mathscr{V}_{\mathscr{R}{C}}=<e_3>$,\\
{\rm (9)} $\alpha-\beta-\gamma= 0,\alpha\neq 0,\beta\neq 0,\gamma\neq 0,\beta\neq\gamma, \lambda=0,\mathscr{V}_{\mathscr{R}{C}}=<e_1>$,\\
{\rm (10)} $-\alpha+\beta-\gamma= 0,\alpha\neq 0,\beta\neq 0,\gamma\neq 0,\alpha\neq\gamma,\lambda=0,  \mathscr{V}_{\mathscr{R}{C}}=<e_2>$,\\
{\rm (11)} $\gamma-\alpha-\beta= 0,\alpha\neq 0,\beta\neq 0,\gamma\neq 0,\alpha\neq\beta, \lambda=0, \mathscr{V}_{\mathscr{R}{C}}=<e_3>$,\\

where $\mathscr{V}_{\mathscr{R}{C}}$ is the vector space of conformal Ricci collineations on $(G_3, g, V)$.
\end{thm}
\begin{Proof}
We know that $\lambda=0$. By the fourth equation, we have $(\beta-\gamma)(\alpha-\beta-\gamma)(\alpha+\beta+\gamma)\lambda_1=0$. \\
case {\rm 1)} If $\alpha+\beta+\gamma= 0$, {\rm (2.18)} holds. We get {\rm (1)}.\\
case {\rm 2)} If $\alpha+\beta+\gamma\neq 0$, $\beta-\gamma=0$, by {\rm (2.18)},
\begin{eqnarray}
\begin{cases}
\alpha(\alpha-\gamma)\lambda_3=0\\[2pt]
\alpha(\alpha-\gamma)\lambda_2=0\\[2pt]
\end{cases}
\end{eqnarray}
case {\rm 2-1)} If $\alpha= 0$, then $\beta=\gamma\neq 0$, {\rm (2.19)} holds. We get {\rm (3)}.\\
case {\rm 2-2)} If $\alpha\neq 0$, then $\alpha=\gamma$. We get {\rm (2)}.\\
case {\rm 2-3)} If $\alpha\neq 0$, $\alpha\neq\gamma$, then $\lambda_2=\lambda_3=0$. We get {\rm (6)}.\\
case {\rm 3)} If $\alpha+\beta+\gamma\neq 0$, $\beta-\gamma\neq 0$, then $\alpha-\beta-\gamma=0$. By {\rm (2.18)},
\begin{eqnarray}
\begin{cases}
(\alpha-\beta)\beta\lambda_3=0\\[2pt]
(\alpha-\gamma)\gamma\lambda_2=0\\[2pt]
\end{cases}
\end{eqnarray}
case {\rm 3-1)} If $\beta= 0$, then $\alpha=\gamma\neq 0$. {\rm (2.20)} holds. We get {\rm (4)}.\\
case {\rm 3-2)} If $\beta\neq 0$, $\alpha-\beta =0$, then $\gamma=0$. {\rm (2.20)} holds. We get {\rm (5)}.\\
case {\rm 3-3)} If $\beta\neq 0$, $\alpha-\beta\neq 0$, we have $\alpha\neq 0,\beta\neq 0,\gamma\neq 0,$ then $\lambda_2=\lambda_3=0$. We get {\rm (9)}.\\
case {\rm 4)} If $\alpha+\beta+\gamma\neq 0$, $\beta-\gamma\neq 0$, $\alpha-\beta-\gamma\neq 0$, then $\lambda_1=0$. By {\rm (2.18)},
\begin{eqnarray}
\begin{cases}
(\alpha-\beta)(\gamma-\alpha-\beta)\lambda_3=0\\[2pt]
(\alpha-\gamma)(-\alpha+\beta-\gamma)\lambda_2=0\\[2pt]
\end{cases}
\end{eqnarray}
case {\rm 4-1)} If $\alpha=\beta$, then $\gamma\neq 0$. By {\rm (2.21)} we have $(\alpha-\gamma)\gamma\lambda_2=0,$ then $\lambda_2=0$. We get {\rm (8)}.\\
case {\rm 4-2)} If $\alpha\neq\beta$, $\gamma-\alpha-\beta =0$, then $\alpha\neq 0,\gamma\neq 0,\beta\neq 0$ . By {\rm (2.)} we have $(\alpha-\gamma)\beta\lambda_2=0$, then $\lambda_2=0$. We get {\rm (11)}.\\
case {\rm 4-3)} If $\alpha\neq\beta$, $\gamma-\alpha-\beta\neq 0$, then $\lambda_3=0$.
If $\alpha=\gamma$, then $\beta\neq 0$. {\rm (2.21)} holds. We get {\rm (7)}.\\
If $\alpha\neq\gamma$, then $-\alpha+\beta-\gamma= 0$, we have $\alpha\neq 0,\beta\neq 0,\gamma\neq 0$. We get {\rm (10)}.\

\end{Proof}

\vskip 0.5 true cm
\noindent{\bf 2.4 Conformal Ricci collineations of $G_4$ }\\
\vskip 0.5 true cm
 By \cite{W}, we have the following Lie algebra of $G_4$ satisfies
\begin{equation}
[e_1,e_2]=-e_2+(2\eta-\beta)e_3,~~\eta=1~{\rm or}-1,~~[e_1,e_3]=-\beta e_2+ e_3,~~[e_2,e_3]=\alpha e_1.
\end{equation}
where $e_1,e_2,e_3$ is a pseudo-orthonormal basis, with $e_3$ timelike.
\begin{lem}The Ricci tensor of  $(G_4,g)$ associated to the Levi-Civita connection $\nabla^{L}$ is determined by
\begin{align}
&Ric^{L}(e_1,e_1)=-\dfrac{\alpha^2}{2},~~~Ric^{L}(e_1,e_2)=0,~~~Ric^{L}(e_1,e_3)=0,\\\notag
&Ric^{L}(e_2,e_2)=\dfrac{\alpha^2}{2}+2\eta(\alpha-\beta)-\alpha\beta+2,~~~Ric^{L}(e_2,e_3)=\alpha+2\eta-\beta,\\\notag
&Ric^{L}(e_3,e_3)=-\dfrac{\alpha^2}{2}-2\beta\eta+\alpha\beta+2.
\end{align}
\end{lem}
Similar to the previous lemma, we can get the Lie derivative of the Ricci tensor $Ric^{L}$ associated to $V$.
\begin{lem}
\begin{align}
&(L_VRic^{L})(e_1,e_1)=0,\\\notag
&(L_VRic^{L})(e_1,e_2)=(-\dfrac{\alpha^2}{2}+\beta^2-2\beta\eta+2)\lambda_2\\\notag
&\hspace{3.5cm}+(-\dfrac{\alpha^3}{2}-\dfrac{\alpha^2\beta}{2}+\alpha\beta^2-2\alpha\beta\eta+2\beta^2\eta+\alpha-3\beta+2\eta)\lambda_3,\\\notag
&(L_VRic^{L})(e_1,e_3)=(\dfrac{\alpha^3}{2}+\dfrac{\alpha^2\beta}{2}-\alpha\beta^2+2\alpha\beta\eta+2\beta^2\eta-\alpha^2\eta-\alpha-5\beta+2\eta)\lambda_2\\\notag
&\hspace{3.5cm}+(-\dfrac{\alpha^2}{2}+\beta^2-4\beta\eta+2)\lambda_3,\\\notag
&(L_VRic^{L})(e_2,e_2)=(\alpha^2-2\beta^2+4\beta\eta+4)\lambda_1,\\\notag
&(L_VRic^{L})(e_2,e_3)=\eta(\alpha-2\beta+2\eta)(\alpha+2\beta-2\eta)\lambda_1,\\\notag
&(L_VRic^{L})(e_3,e_3)=(\alpha^2-2\beta^2+8\beta\eta-4)\lambda_1.
\end{align}
\end{lem}
Then, if $V$ is a conformal Ricci collineation associated to the Levi-Civita connection, by Lemma 2.12 and Theorem 2.1, we have the following equations:
\begin{eqnarray}
       \begin{cases}
       0=2\lambda\\[2pt]
       (-\dfrac{\alpha^2}{2}+\beta^2-2\beta\eta+2)\lambda_2\\[2pt]
       \hspace{1cm}+(-\dfrac{\alpha^3}{2}-\dfrac{\alpha^2\beta}{2}+\alpha\beta^2-2\alpha\beta\eta+2\beta^2\eta+\alpha-3\beta+2\eta)\lambda_3=0\\[2pt]
       (\dfrac{\alpha^3}{2}+\dfrac{\alpha^2\beta}{2}-\alpha\beta^2+2\alpha\beta\eta+2\beta^2\eta-\alpha^2\eta-\alpha-5\beta+2\eta)\lambda_2\\[2pt]
       \hspace{1cm}+(-\dfrac{\alpha^2}{2}+\beta^2-4\beta\eta+2)\lambda_3=0\\[2pt]
       (\alpha^2-2\beta^2+4\beta\eta+4)\lambda_1=2\lambda\\[2pt]
       \eta(\alpha-2\beta+2\eta)(\alpha+2\beta-2\eta)\lambda_1=0\\[2pt]
       (\alpha^2-2\beta^2+8\beta\eta-4)\lambda_1=-2\lambda\\[2pt]
       \end{cases}
\end{eqnarray}
By solving (2.25) , we get
\begin{thm}
$(G_4, g, V)$ admits conformal Ricci collineations associated to the Levi-Civita connection if and only if one of the following holds:\\
{\rm (1)} $\eta=1~{\rm or}~-1$, $\lambda=0$, $\alpha=2\eta$, $\beta=0$, $\mathscr{V}_{\mathscr{R}{C}}=<e_2,e_3>$,\\
{\rm (2)} $\eta=1~{\rm or}~-1$, $\lambda=0$, $m= 0$, $n= 0$, $\alpha^2-4\neq 0$, $\mathscr{V}_{\mathscr{R}{C}}=<\dfrac{2\beta\eta}{\alpha^2-4}e_2+e_3>$,\\
{\rm (3)} $\eta=1~{\rm or}~-1$, $\lambda=0$, $m= 0$, $n\neq 0$, $\alpha^2-4-n=0$, $\mathscr{V}_{\mathscr{R}{C}}=<e_2>$.\\
{\rm (4)} $\eta=1~{\rm or}~-1$, $\lambda=0$, $m\neq 0$, $ m(m-2\beta)-n(4\beta^2-\alpha^2-8\beta+4-n)=0$,\\ $\mathscr{V}_{\mathscr{R}{C}}=<-\dfrac{n}{m}e_2+e_3>$.\\
where $m=-\dfrac{\alpha^2}{2}+\beta^2-2\beta\eta+2, n=-\dfrac{\alpha^3}{2}-\dfrac{\alpha^2\beta}{2}+\alpha\beta^2-2\alpha\beta\eta+2\beta^2\eta+\alpha-3\beta+2\eta$, $\mathscr{V}_{\mathscr{R}{C}}$ is the vector space of conformal Ricci collineations on $(G_4, g, V)$.
\end{thm}
\begin{Proof}
We know that $\eta=1~{\rm or}~-1$. By the fifth equation, we have $(\alpha-2\beta+2\eta)(\alpha+2\beta-2\eta)\lambda_1=0$. \\
Then by the sixth equation, $2\beta^2\lambda_1=0$. If $\beta=0$, then $(\alpha^2+4)\lambda_1=0$, i.e. $\lambda_1=0$; otherwise $\beta\neq 0$, then $\lambda_1=0$. So $\lambda=\lambda_1=0$,  let $m=-\dfrac{\alpha^2}{2}+\beta^2-2\beta\eta+2, n=-\dfrac{\alpha^3}{2}-\dfrac{\alpha^2\beta}{2}+\alpha\beta^2-2\alpha\beta\eta+2\beta^2\eta+\alpha-3\beta+2\eta$, {\rm (2.25)} can be simplified to
\begin{eqnarray}
\begin{cases}
       m\lambda_2+n\lambda_3=0\\[2pt]
       (4\beta^2\eta-\alpha^2\eta-8\beta+4\eta-n)\lambda_2+(m-2\beta\eta)\lambda_3=0\\[2pt]
\end{cases}
\end{eqnarray}
case {\rm 1)} If $\eta=1$, let $m_1=m$, $n_1=n$,  by {\rm (2.)} we have
\begin{eqnarray}
\begin{cases}
       m_1\lambda_2+n_1\lambda_3=0\\[2pt]
       (4\beta^2-\alpha^2-8\beta+4-n_1)\lambda_2+(m_1-2\beta)\lambda_3=0\\[2pt]
\end{cases}
\end{eqnarray}
case {\rm 1-1)} If $m_1= 0$, then $n_1\lambda_3=0$.\\
case {\rm 1-1-1)} If $n_1= 0$, by $m_1=-\dfrac{\alpha^2}{2}+(\beta-1)^2+1=0$, {\rm (2.27)} can be simplified to
\begin{align}\notag
&(-\alpha^2+4(\beta-1)^2)\lambda_2-2\beta\lambda_3=(\alpha^2-4)\lambda_2-2\beta\lambda_3=0.
\end{align}
If $\alpha^2-4=0$, then $\beta=0~{\rm or}~2$, $ 2\beta\lambda_3=0$. When $\beta=0$, by $n_1= 0$ we have $\alpha=2$, this situation falls into {\rm (1)}; otherwise $\beta=2$, by $n_1= 0$, we get a contradiction.\\
If $\alpha^2-4\neq 0$,  then $\lambda_2=\dfrac{2\beta}{\alpha^2-4}\lambda_3$. This situation falls into {\rm (2)}.\\
case {\rm 1-1-2)} If $n_1\neq 0$, then $\lambda_3=0$. By {\rm (2.27)} we have
\begin{align}\notag
&(\alpha^2-4-n_1)\lambda_2=0,
\end{align}
then $\alpha^2-4-n_1=0.$ This situation falls into {\rm (3)}.\\
case {\rm 1-2)} If $m_1\neq 0$,  then $\lambda_2=-\dfrac{n_1}{m_1}\lambda_3$, $ m_1(m_1-2\beta)-n_1(4\beta^2-\alpha^2-8\beta+4-n_1)=0$. This situation falls into {\rm (4)}.\\
case {\rm 2)} If $\eta=-1$, similar to $\eta=1$, $\lambda=\lambda_1=0$, let $m_2=m$, $n_2=n$, we have\\
case {\rm 2-1)} $\alpha=-2$, $\beta=0$, this situation falls into {\rm (1)}.\\
case {\rm 2-2)} $m_2=0$, $n_2=0$, $\alpha^2-4\neq 0$, $\lambda_2=-\dfrac{2\beta}{\alpha^2-4}\lambda_3$. This situation falls into {\rm (2)}.\\
case {\rm 2-3)} $m_2= 0$, $n_2\neq 0$, $\alpha^2-4-n_2=0$, $\lambda_3=0$. This situation falls into {\rm (3)}.\\
case {\rm 2-4)} $m_2\neq 0$, $ m_2(m_2-2\beta)-n_2(4\beta^2-\alpha^2-8\beta+4-n_2)=0$, $\lambda_2=-\dfrac{n_2}{m_2}\lambda_3$. This situation falls into {\rm (4)}.
\end{Proof}

\section{  Conformal Ricci collineations associated to the Levi-Civita connection on three-dimensional Lorentzian non-unimodular Lie groups}
\vskip 0.5 true cm
\noindent{\bf 3.1 Conformal Ricci collineations of $G_5$ }\\
\vskip 0.5 true cm
 By \cite{W}, we have the following Lie algebra of $G_5$ satisfies
\begin{equation}
[e_1,e_2]=0,~~[e_1,e_3]=\alpha e_1+\beta e_2,~~[e_2,e_3]=\gamma e_1+\delta e_2,~~\alpha+\delta\neq 0,~~\alpha\gamma+\beta\delta=0.
\end{equation}
where $e_1,e_2,e_3$ is a pseudo-orthonormal basis, with $e_3$ timelike.
\begin{lem}The Ricci tensor of  $(G_5,g)$ associated to the Levi-Civita connection $\nabla^{L}$ is determined by
\begin{align}
&Ric^{L}(e_1,e_1)=\alpha^2+\alpha\delta+\dfrac{\beta^2-\gamma^2}{2},~~~Ric^{L}(e_1,e_2)=0,~~~Ric^{L}(e_1,e_3)=0,\\\notag
&Ric^{L}(e_2,e_2)=\delta^2+\alpha\delta-\dfrac{\beta^2-\gamma^2}{2},~~~Ric^{L}(e_2,e_3)=0,~~~Ric^{L}(e_3,e_3)=-(\alpha^2+\delta^2+\dfrac{(\beta+\gamma)^2}{2}).
\end{align}
\end{lem}
Similar to the previous lemma, we can get the Lie derivative of the Ricci tensor $Ric^{L}$ associated to $V$.
\begin{lem}
\begin{align}
&(L_VRic^{L})(e_1,e_1)=2\alpha(\alpha^2+\alpha\delta+\dfrac{\beta^2-\gamma^2}{2})\lambda_3,\\\notag
&(L_VRic^{L})(e_1,e_2)=[\beta(\delta^2+\alpha\delta-\dfrac{\beta^2-\gamma^2}{2})+\gamma(\alpha^2+\alpha\delta+\dfrac{\beta^2-\gamma^2}{2})]\lambda_3,\\\notag
&(L_VRic^{L})(e_1,e_3)=-\alpha(\alpha^2+\alpha\delta+\dfrac{\beta^2-\gamma^2}{2})\lambda_1-\gamma(\alpha^2+\alpha\delta+\dfrac{\beta^2-\gamma^2}{2})\lambda_2,\\\notag
&(L_VRic^{L})(e_2,e_2)=2\delta(\delta^2+\alpha\delta-\dfrac{\beta^2-\gamma^2}{2})\lambda_3,\\\notag
&(L_VRic^{L})(e_2,e_3)=-\beta(\delta^2+\alpha\delta-\dfrac{\beta^2-\gamma^2}{2})\lambda_1-\delta(\delta^2+\alpha\delta-\dfrac{\beta^2-\gamma^2}{2})\lambda_2,\\\notag
&(L_VRic^{L})(e_3,e_3)=0.
\end{align}
\end{lem}
Then, if $V$ is a conformal Ricci collineation associated to the Levi-Civita connection, by Lemma 3.2 and Theorem 2.1, we have the following equations:
\begin{eqnarray}
       \begin{cases}
       2\alpha(\alpha^2+\alpha\delta+\dfrac{\beta^2-\gamma^2}{2})\lambda_3=2\lambda\\[2pt]
       [\beta(\delta^2+\alpha\delta-\dfrac{\beta^2-\gamma^2}{2})+\gamma(\alpha^2+\alpha\delta+\dfrac{\beta^2-\gamma^2}{2})]\lambda_3=0\\[2pt]
       (\alpha^2+\alpha\delta+\dfrac{\beta^2-\gamma^2}{2})(\alpha\lambda_1+\gamma\lambda_2)=0\\[2pt]
       2\delta(\delta^2+\alpha\delta-\dfrac{\beta^2-\gamma^2}{2})\lambda_3=2\lambda\\[2pt]
       (\delta^2+\alpha\delta-\dfrac{\beta^2-\gamma^2}{2})(\beta\lambda_1+\delta\lambda_2)=0\\[2pt]
       0=-2\lambda\\[2pt]
       \end{cases}
\end{eqnarray}
By solving (3.4) , we get
\begin{thm}
$(G_5, g, V)$ admits conformal Ricci collineations associated to the Levi-Civita connection if and only if one of the following holds:\\
{\rm (1)} $\alpha^2+\alpha\delta-\dfrac{\gamma^2}{2}= 0,\beta= 0,\alpha+\delta\neq 0,\alpha\gamma=0,\delta\neq 0,\lambda=0,~~ \mathscr{V}_{\mathscr{R}{C}}=<e_1>$,\\
{\rm (2)} $\alpha^2+\alpha\delta+\dfrac{\beta^2-\gamma^2}{2}= 0,\beta\neq 0,\alpha+\delta\neq 0,\alpha\gamma+\beta\delta=0,\lambda=0,~~ \mathscr{V}_{\mathscr{R}{C}}=<-\dfrac{\delta}{\beta}e_1+e_2>$,\\
{\rm (3)} $\alpha= 0,\beta= 0,\delta\neq 0,\gamma\neq 0,\lambda=0,~~ \mathscr{V}_{\mathscr{R}{C}}=<e_1>$,\\
{\rm (4)} $\delta^2+\alpha\delta-\dfrac{\beta^2-\gamma^2}{2}= 0,\alpha\neq 0,\alpha+\delta\neq 0,\alpha\gamma+\beta\delta=0,\lambda=0,~~ \mathscr{V}_{\mathscr{R}{C}}=<-\dfrac{\gamma}{\alpha}e_1+e_2>$,\\
{\rm (5)} $\alpha^2+\alpha\delta+\dfrac{\beta^2-\gamma^2}{2}\neq 0,\delta^2+\alpha\delta-\dfrac{\beta^2-\gamma^2}{2}\neq 0,\alpha\neq 0,\alpha+\delta\neq 0,\alpha\gamma+\beta\delta=0,\alpha\delta-\beta\gamma=0,\lambda=0,~~ \mathscr{V}_{\mathscr{R}{C}}=<-\dfrac{\gamma}{\alpha}e_1+e_2>$.\\
where $\mathscr{V}_{\mathscr{R}{C}}$ is the vector space of conformal Ricci collineations on $(G_5, g, V)$.
\end{thm}
\begin{Proof}
We know that $\lambda= 0,~\alpha+\delta\neq 0,~\alpha\gamma+\beta\delta=0$. By the fourth equation, we have $\delta(\delta^2+\alpha\delta-\dfrac{\beta^2-\gamma^2}{2})\lambda_3=0$.\\
case {\rm 1)} If $\lambda_3=0$, {\rm (3.4)} can be simplified to
\begin{eqnarray}
\begin{cases}
      (\alpha^2+\alpha\delta+\dfrac{\beta^2-\gamma^2}{2})(\alpha\lambda_1+\gamma\lambda_2)=0\\[2pt]
       (\delta^2+\alpha\delta-\dfrac{\beta^2-\gamma^2}{2})(\beta\lambda_1+\delta\lambda_2)=0\\[2pt]
\end{cases}
\end{eqnarray}
case {\rm 1-1)} If $\alpha^2+\alpha\delta+\dfrac{\beta^2-\gamma^2}{2}= 0$, {\rm (3.5)} can be simplified to
\begin{align}\notag
&(\alpha+\delta)^2(\beta\lambda_1+\delta\lambda_2)=0,
\end{align}
i.e. $\beta\lambda_1+\delta\lambda_2=0$.\\
case {\rm 1-1-1)} If $\beta=0$, then $\delta\lambda_2=0$. When $\delta=0$, by $\alpha+\delta=0,~\alpha\gamma+\beta\delta=0$, we have $\alpha\neq 0,~\gamma=0,~\alpha^2+\alpha\delta+\dfrac{\beta^2-\gamma^2}{2}=\alpha^2\neq 0$, we get a contradiction; otherwise $\delta\neq 0$, $\lambda_2=0$. we get {\rm (1)}.\\
case {\rm 1-1-2)} If $\beta\neq 0$,  then $\lambda_1=-\dfrac{\delta}{\beta}\lambda_2$. We get {\rm (2)}.\\
case {\rm 1-2)} If $\alpha^2+\alpha\delta+\dfrac{\beta^2-\gamma^2}{2}\neq 0$, then $\alpha\lambda_1+\gamma\lambda_2=0$.\\
case {\rm 1-2-1)} If $\alpha=0$, by $\alpha+\delta=0,~\alpha\gamma+\beta\delta=0$, we have $\delta\neq 0,~\beta=0,~\gamma\neq 0$, then $\lambda_2=0$. We get {\rm (3)}.\\
case {\rm 1-2-2)} If $\alpha\neq 0$,  then $\lambda_1=-\dfrac{\gamma}{\alpha}\lambda_2$. {\rm (3.5)} can be simplified to
\begin{align}\notag
&(\delta^2+\alpha\delta+\dfrac{\beta^2-\gamma^2}{2})(-\beta\gamma+\alpha\delta)\lambda_2=0.
\end{align}
If $\delta^2+\alpha\delta+\dfrac{\beta^2-\gamma^2}{2}=0$, {\rm (3.5)} holds. we get {\rm (4)}.\\
If $\delta^2+\alpha\delta+\dfrac{\beta^2-\gamma^2}{2}\neq 0$, then $\alpha\delta-\beta\gamma=0$. We get {\rm (5)}.\\
case {\rm 2)} If $\lambda_3\neq 0$, $\delta=0$, we have $\alpha\neq 0,~\gamma=0$, by {\rm (3.4)} we get $\beta=0,~\alpha=0$. This is a contradiction.\\
So $\delta\neq 0$, $\delta^2+\alpha\delta-\dfrac{\beta^2-\gamma^2}{2}= 0$, then $\alpha^2+\alpha\delta+\dfrac{\beta^2-\gamma^2}{2}\neq 0$. By {\rm (3.4)} we have $\alpha=\gamma=0$, by $\alpha\gamma-\beta\delta=0$, we have $\beta=0,$ so$\delta^2+\alpha\delta+\dfrac{\beta^2-\gamma^2}{2}\neq 0$, we get a contradiction.
\end{Proof}

\vskip 0.5 true cm
\noindent{\bf 3.2 Conformal Ricci collineations of $G_6$ }\\
\vskip 0.5 true cm
 By \cite{W}, we have the following Lie algebra of $G_6$ satisfies
\begin{equation}
[e_1,e_2]=\alpha e_2+\beta e_3,~~[e_1,e_3]=\gamma e_2+\delta e_3,~~[e_2,e_3]=0,~~\alpha+\delta\neq 0,~~\alpha\gamma-\beta\delta=0.
\end{equation}
where $e_1,e_2,e_3$ is a pseudo-orthonormal basis, with $e_3$ timelike.
\begin{lem}The Ricci tensor of  $(G_6,g)$ associated to the Levi-Civita connection $\nabla^{L}$ is determined by
\begin{align}
&Ric^{L}(e_1,e_1)=-\alpha^2-\delta^2+\dfrac{(\beta-\gamma)^2}{2},~~~Ric^{L}(e_1,e_2)=0,~~~Ric^{L}(e_1,e_3)=0,\\\notag
&Ric^{L}(e_2,e_2)=-\alpha^2-\alpha\delta+\dfrac{\beta^2-\gamma^2}{2},~~~Ric^{L}(e_2,e_3)=0,~~~Ric^{L}(e_3,e_3)=\delta^2+\alpha\delta+\dfrac{\beta^2-\gamma^2}{2}.
\end{align}
\end{lem}
Similar to the previous lemma, we can get the Lie derivative of the Ricci tensor $Ric^{L}$ associated to $V$.
\begin{lem}
\begin{align}
&(L_VRic^{L})(e_1,e_1)=0,\\\notag
&(L_VRic^{L})(e_1,e_2)=\alpha(-\alpha^2-\alpha\delta+\dfrac{\beta^2-\gamma^2}{2})\lambda_2+\gamma(-\alpha^2-\alpha\delta+\dfrac{\beta^2-\gamma^2}{2})\lambda_3,\\\notag
&(L_VRic^{L})(e_1,e_3)=\beta(\delta^2+\alpha\delta+\dfrac{\beta^2-\gamma^2}{2})\lambda_2+\delta(\delta^2+\alpha\delta+\dfrac{\beta^2-\gamma^2}{2})\lambda_3,\\\notag
&(L_VRic^{L})(e_2,e_2)=-2\alpha(-\alpha^2-\alpha\delta+\dfrac{\beta^2-\gamma^2}{2})\lambda_1,\\\notag
&(L_VRic^{L})(e_2,e_3)=[-\beta(\delta^2+\alpha\delta+\dfrac{\beta^2-\gamma^2}{2})-\gamma(-\alpha^2-\alpha\delta+\dfrac{\beta^2-\gamma^2}{2})]\lambda_1,\\\notag
&(L_VRic^{L})(e_3,e_3)=-2\delta(\delta^2+\alpha\delta+\dfrac{\beta^2-\gamma^2}{2})\lambda_1.
\end{align}
\end{lem}
Then, if $V$ is a conformal Ricci collineation associated to the Levi-Civita connection, by Lemma 3.5 and Theorem 2.1, we have the following equations:
\begin{eqnarray}
       \begin{cases}
       0=2\lambda\\[2pt]
       (-\alpha^2-\alpha\delta+\dfrac{\beta^2-\gamma^2}{2})(\alpha\lambda_2+\gamma\lambda_3)=0\\[2pt]
       (\delta^2+\alpha\delta+\dfrac{\beta^2-\gamma^2}{2})(\beta\lambda_2+\delta\lambda_3)=0\\[2pt]
       -2\alpha(-\alpha^2-\alpha\delta+\dfrac{\beta^2-\gamma^2}{2})\lambda_1=2\lambda\\[2pt]
       [-\beta(\delta^2+\alpha\delta+\dfrac{\beta^2-\gamma^2}{2})-\gamma(-\alpha^2-\alpha\delta+\dfrac{\beta^2-\gamma^2}{2})]\lambda_1=0\\[2pt]
       -2\delta(\delta^2+\alpha\delta+\dfrac{\beta^2-\gamma^2}{2})\lambda_1=-2\lambda\\[2pt]
       \end{cases}
\end{eqnarray}
By solving (3.9) , we get
\begin{thm}
$(G_6, g, V)$ admits conformal Ricci collineations associated to the Levi-Civita connection if and only if one of the following holds:\\
{\rm (1)} $\alpha^2+\alpha\delta+\dfrac{\gamma^2}{2}= 0,\beta= 0,\alpha+\delta\neq 0,\alpha\gamma=0,\delta\neq 0,\lambda=0,~~ \mathscr{V}_{\mathscr{R}{C}}=<e_2>$,\\
{\rm (2)} $\alpha^2+\alpha\delta-\dfrac{\beta^2-\gamma^2}{2}= 0,\beta\neq 0,\alpha+\delta\neq 0,\alpha\gamma-\beta\delta=0,\lambda=0,~~ \mathscr{V}_{\mathscr{R}{C}}=<-\dfrac{\delta}{\beta}e_2+e_3>$,\\
{\rm (3)} $\alpha= 0,\beta= 0,\delta\neq 0,\gamma\neq 0,\lambda=0,~~ \mathscr{V}_{\mathscr{R}{C}}=<e_2>$,\\
{\rm (4)} $\delta^2+\alpha\delta+\dfrac{\beta^2-\gamma^2}{2}= 0,\alpha\neq 0,\alpha+\delta\neq 0,\alpha\gamma-\beta\delta=0,\lambda=0,~~ \mathscr{V}_{\mathscr{R}{C}}=<-\dfrac{\gamma}{\alpha}e_2+e_3>$,\\
{\rm (5)} $\alpha^2+\alpha\delta-\dfrac{\beta^2-\gamma^2}{2}\neq 0,\delta^2+\alpha\delta+\dfrac{\beta^2-\gamma^2}{2}\neq 0,\alpha\neq 0,\alpha+\delta\neq 0,\alpha\gamma-\beta\delta=0,\alpha\delta-\beta\gamma=0,\lambda=0,~~ \mathscr{V}_{\mathscr{R}{C}}=<-\dfrac{\gamma}{\alpha}e_2+e_3>$.\\
where $\mathscr{V}_{\mathscr{R}{C}}$ is the vector space of conformal Ricci collineations on $(G_6, g, V)$.
\end{thm}
\begin{Proof}
We know that $\lambda= 0,~\alpha+\delta\neq 0,~\alpha\gamma-\beta\delta=0$. By the fifth equation, we have $\alpha(\alpha^2+\alpha\delta-\dfrac{\beta^2-\gamma^2}{2})\lambda_1=0$.\\
case {\rm 1)} If $\lambda_1=0$, {\rm (3.9)} can be simplified to
\begin{eqnarray}
\begin{cases}
      (\alpha^2+\alpha\delta-\dfrac{\beta^2-\gamma^2}{2})(\alpha\lambda_2+\gamma\lambda_3)=0\\[2pt]
       (\delta^2+\alpha\delta+\dfrac{\beta^2-\gamma^2}{2})(\beta\lambda_2+\delta\lambda_3)=0\\[2pt]
\end{cases}
\end{eqnarray}
case {\rm 1-1)} If $\alpha^2+\alpha\delta-\dfrac{\beta^2-\gamma^2}{2}= 0$, {\rm (3.10)} can be simplified to
\begin{align}\notag
&(\alpha+\delta)^2(\beta\lambda_2+\delta\lambda_3)=0,
\end{align}
i.e. $\beta\lambda_2+\delta\lambda_3=0$.\\
case {\rm 1-1-1)} If $\beta=0$, then $\delta\lambda_3=0$. When $\delta=0$, by $\alpha+\delta=0,~\alpha\gamma-\beta\delta=0$, we have $\alpha\neq 0,~\gamma=0,~\alpha^2+\alpha\delta-\dfrac{\beta^2-\gamma^2}{2}\neq 0$, we get a contradiction; otherwise $\delta\neq 0$, $\lambda_3=0$. We get {\rm (1)}.\\
case {\rm 1-1-2)} If $\beta\neq 0$,  then $\lambda_2=-\dfrac{\delta}{\beta}\lambda_3$. We get {\rm (2)}.\\
case {\rm 1-2)} If $\alpha^2+\alpha\delta-\dfrac{\beta^2-\gamma^2}{2}\neq 0$, then $\alpha\lambda_2+\gamma\lambda_3=0$.\\
case {\rm 1-2-1)} If $\alpha=0$, by $\alpha+\delta\neq 0,~\alpha\gamma-\beta\delta=0$, we have $\delta\neq 0,~\beta=0,~\gamma\neq 0$, then $\lambda_3=0$. We get {\rm (3)}.\\
case {\rm 1-2-2)} If $\alpha\neq 0$,  then $\lambda_2=-\dfrac{\gamma}{\alpha}\lambda_3$. {\rm (3.10)} can be simplified to
\begin{align}\notag
&(\delta^2+\alpha\delta+\dfrac{\beta^2-\gamma^2}{2})(-\beta\gamma+\alpha\delta)\lambda_3=0.
\end{align}
If $\delta^2+\alpha\delta+\dfrac{\beta^2-\gamma^2}{2}=0$, {\rm (3.10)} holds. We get {\rm (4)}.\\
If $\delta^2+\alpha\delta+\dfrac{\beta^2-\gamma^2}{2}\neq 0$,  then $\alpha\delta-\beta\gamma=0$. We get {\rm (5)}.\\
case {\rm 2)} If $\lambda_1\neq 0$, $\alpha=0$, we have $\delta\neq 0,~\beta=0$, by {\rm (3.9)} we get $\gamma=0,~\delta=0$. This is a contradiction.\\
So $\alpha\neq 0$, $\alpha^2+\alpha\delta-\dfrac{\beta^2-\gamma^2}{2}= 0$, then $\delta^2+\alpha\delta+\dfrac{\beta^2-\gamma^2}{2}\neq 0$.  By {\rm (3.9)} we have $\beta=\delta=0$, by $\alpha\gamma-\beta\delta=0$, we have $\gamma=0,$ so$\alpha^2+\alpha\delta-\dfrac{\beta^2-\gamma^2}{2}=\alpha^2\neq 0$,, we get a contradiction.
\end{Proof}

\vskip 0.5 true cm
\noindent{\bf 3.3 Conformal Ricci collineations of $G_7$ }\\
\vskip 0.5 true cm
 By \cite{W}, we have the following Lie algebra of $G_7$ satisfies
\begin{equation}
[e_1,e_2]=-\alpha e_1-\beta e_2-\beta e_3,~~[e_1,e_3]=\alpha e_1+\beta e_2+\beta e_3,~~[e_2,e_3]=\gamma e_1+\delta e_2+\delta e_3,~~\alpha+\delta\neq 0,~~\alpha\gamma=0.
\end{equation}
where $e_1,e_2,e_3$ is a pseudo-orthonormal basis, with $e_3$ timelike.
\begin{lem}The Ricci tensor of  $(G_7,g)$ associated to the Levi-Civita connection $\nabla^{L}$ is determined by
\begin{align}
&Ric^{L}(e_1,e_1)=-\dfrac{\gamma^2}{2},~~~Ric^{L}(e_1,e_2)=0,~~~Ric^{L}(e_1,e_3)=0,\\\notag
&Ric^{L}(e_2,e_2)=-\alpha^2-\dfrac{\gamma^2}{2}+\alpha\delta-\beta\gamma,~~~Ric^{L}(e_2,e_3)=\alpha^2-\alpha\delta+\beta\gamma,~~~Ric^{L}(e_3,e_3)=-\alpha^2-\dfrac{\gamma^2}{2}+\alpha\delta-\beta\gamma.
\end{align}
\end{lem}
Similar to the previous lemma, we can get the Lie derivative of the Ricci tensor $Ric^{L}$ associated to $V$.
\begin{lem}
\begin{align}
&(L_VRic^{L})(e_1,e_1)=0,\\\notag
&(L_VRic^{L})(e_1,e_2)=-\dfrac{\beta\gamma^2}{2}\lambda_2+(\dfrac{\beta\gamma^2}{2}-\dfrac{\gamma^3}{2})\lambda_3,\\\notag
&(L_VRic^{L})(e_1,e_3)=(\dfrac{\beta\gamma^2}{2}+\dfrac{\gamma^3}{2})\lambda_2-\dfrac{\beta\gamma^2}{2}\lambda_3,\\\notag
&(L_VRic^{L})(e_2,e_2)=\beta\gamma^2\lambda_1+\delta\gamma^2\lambda_3,\\\notag
&(L_VRic^{L})(e_2,e_3)=-\beta\gamma^2\lambda_1-\dfrac{\delta\gamma^2}{2}\lambda_2-\dfrac{\delta\gamma^2}{2}\lambda_3,\\\notag
&(L_VRic^{L})(e_3,e_3)=\beta\gamma^2\lambda_1+\delta\gamma^2\lambda_2.
\end{align}
\end{lem}
Then, if $V$ is a conformal Ricci collineation associated to the Levi-Civita connection, by Lemma 3.8 and Theorem 2.1, we have the following equations:
\begin{eqnarray}
       \begin{cases}
       0=2\lambda\\[2pt]
       -\dfrac{\beta\gamma^2}{2}\lambda_2+(\dfrac{\beta\gamma^2}{2}-\dfrac{\gamma^3}{2})\lambda_3=0\\[2pt]
       (\dfrac{\beta\gamma^2}{2}+\dfrac{\gamma^3}{2})\lambda_2-\dfrac{\beta\gamma^2}{2}\lambda_3=0\\[2pt]
       \beta\gamma^2\lambda_1+\delta\gamma^2\lambda_3=2\lambda\\[2pt]
       -\beta\gamma^2\lambda_1-\dfrac{\delta\gamma^2}{2}\lambda_2-\dfrac{\delta\gamma^2}{2}\lambda_3=0\\[2pt]
       \beta\gamma^2\lambda_1+\delta\gamma^2\lambda_2=-2\lambda\\[2pt]
       \end{cases}
\end{eqnarray}
By solving (3.14) , we get
\begin{thm}
$(G_7, g, V)$ admits conformal Ricci collineations associated to the Levi-Civita connection if and only if one of the following holds:\\
{\rm (1)} $\gamma= 0,\alpha+\delta\neq 0,\lambda=0,~ \mathscr{V}_{\mathscr{R}{C}}=<e_1,e_2,e_3>$,\\
{\rm (2)} $\alpha=0,\beta= 0,\delta\neq 0,\gamma\neq 0,\lambda=0, ~\mathscr{V}_{\mathscr{R}{C}}=<e_1>$.\\
where $\mathscr{V}_{\mathscr{R}{C}}$ is the vector space of conformal Ricci collineations on $(G_7, g, V)$.
\end{thm}
\begin{Proof}
We know that $\lambda=0,~\alpha+\delta\neq 0,~\alpha\gamma=0.$ By {\rm (3.14)}, we have $\delta\gamma^2(\lambda_2-\lambda_3)=0$. \\
case {\rm 1)} If $\delta= 0$, by $\alpha+\delta\neq 0,~\alpha\gamma=0$, we have $\alpha\neq 0,~\gamma=0$, {\rm (3.14)} holds. This situation falls into {\rm (1)}.\\
case {\rm 2)} If $\delta\neq 0$, $\gamma=0$, $\alpha+\delta\neq 0$, {\rm (3.14)} holds. This situation falls into {\rm (1)}.\\
case {\rm 3)} If $\delta\neq 0$, $\gamma\neq 0$, then $\lambda_2=\lambda_3$, by $\alpha\gamma=0$, we have $\alpha=0$. {\rm (3.14)} can be simplified to
\begin{eqnarray}
\begin{cases}
      \dfrac{\gamma^3}{2}\lambda_2=0\\[2pt]
      \beta\gamma^2\lambda_1+\delta\gamma^2\lambda_2=0\\[2pt]
\end{cases}
\end{eqnarray}
We get $\lambda_2=\lambda_3=0$, $\beta=0$.  We get {\rm (2)}.\\
\end{Proof}

\section{Acknowledgements}

The author was supported in part by NSFC No.11771070.

\vskip 1 true cm


\bigskip
\bigskip

\noindent {\footnotesize {\it Yanli. Wang} \\
{School of Mathematics and Statistics, Northeast Normal University, Changchun 130024, China}\\
{Email: wangyl553@nenu.edu.cn}

\end{document}